# ON MARTINGALE APPROXIMATIONS[1]

## By Ou Zhao and Michael Woodroofe

### *Yale University and University of Michigan*

Consider additive functionals of a Markov chain $W_k$, with stationary (marginal) distribution and transition function denoted by $\pi$ and $Q$, say $S_n = g(W_1) + \cdots + g(W_n)$, where $g$ is square integrable and has mean 0 with respect to $\pi$. If $S_n$ has the form $S_n = M_n + R_n$, where $M_n$ is a square integrable martingale with stationary increments and $E(R_n^2) = o(n)$, then $g$ is said to admit a martingale approximation. Necessary and sufficient conditions for such an approximation are developed. Two obvious necessary conditions are $E[E(S_n|W_1)^2] = o(n)$ and $\lim_{n\to\infty} E(S_n^2)/n < \infty$. Assuming the first of these, let $\|g\|_+^2 = \limsup_{n\to\infty} E(S_n^2)/n$; then $\|\cdot\|_+$ defines a pseudo norm on the subspace of $L^2(\pi)$ where it is finite. In one main result, a simple necessary and sufficient condition for a martingale approximation is developed in terms of $\|\cdot\|_+$. Let $Q^*$ denote the adjoint operator to $Q$, regarded as a linear operator from $L^2(\pi)$ into itself, and consider co-isometries ($QQ^* = I$), an important special case that includes shift processes. In another main result a convenient orthonormal basis for $L_0^2(\pi)$ is identified along with a simple necessary and sufficient condition for the existence of a martingale approximation in terms of the coefficients of the expansion of $g$ with respect to this basis.

## 1. Introduction.

Some notation is necessary to describe the results of the paper. Let $\ldots, W_{-1}, W_0, W_1, \ldots$ denote a stationary, ergodic Markov chain with values in a measurable space $\mathcal{W}$. The marginal distribution and transition function of the chain are denoted by $\pi$ and $Q$; thus, $\pi\{B\} = P[W_n \in B]$ and $Q(w; B) = P[W_{n+1} \in B | W_n = w]$ for $w \in \mathcal{W}$ and measurable sets $B \subseteq \mathcal{W}$. In addition, $Q$ denotes the operator, defined by

$$Qf(w) = \int_{\mathcal{W}} f(z)Q(w; dz) \qquad \text{a.e. } (\pi)$$

Received August 2007; revised December 2007.

[1]Supported by the National Science Foundation.

for $f \in L^1(\pi)$, and the iterates of $Q$ are denoted by $Q^k = Q \circ \cdots \circ Q$ ($k$ times). Thus, $Q^k f(w) = E[f(W_{n+k})|W_n = w]$ a.e. ($\pi$) for $f \in L^1(\pi)$. The probability space on which $\ldots, W_{-1}, W_0, W_1, \ldots$ are defined is denoted by $(\Omega, \mathcal{A}, P)$, and $\mathcal{F}_n = \sigma\{\ldots, W_{n-1}, W_n\}$. Finally, $\|\cdot\|$ and $\langle \cdot, \cdot \rangle$ denote the norm and inner product in an $L^2$ space, which may vary from one usage to the next.

Observe that no stringent conditions, like Harris recurrence or even irreducibility, have been placed on the Markov chain. In particular, if $\ldots, \xi_{-1}, \xi_0,$ $\xi_1, \ldots$ are i.i.d. with common distribution $\rho$, say, then the *shift process* $W_k = (\ldots, \xi_{k-1}, \xi_k)$ satisfies the conditions placed on the chain with $\pi = \rho^{\mathbb{N}}$, where $\mathbb{N} = \{0, 1, 2, \ldots\}$, and $Qg(w) = \int g(w, x)\rho\{dx\}$ for $g \in L^1(\pi)$. Shift processes abound in books on time series—for example, [2] and [16].

Next let $L_0^2(\pi)$ be the set of $g \in L^2(\pi)$ for which $\int_{\mathcal{W}} g \, d\pi = 0$; and, for $g \in L_0^2(\pi)$, consider stationary sequences of the form $X_k = g(W_k)$ and their sums $S_n = X_1 + \cdots + X_n$. Thus,

$$S_n = S_n(g) = g(W_1) + \cdots + g(W_n).$$

The question addressed here is the existence of a martingale $M_1, M_2, \ldots$ with respect to $\mathcal{F}_0, \mathcal{F}_1, \mathcal{F}_2, \ldots$ having stationary increments and a sequence of remainder terms $R_1, R_2, \ldots$ for which $\|R_n\| = o(\sqrt{n})$ and

$$(1) \qquad\qquad\qquad S_n = M_n + R_n.$$

If (1) holds, we say that $g$ admits a *martingale approximation*. Ever since the work of Gordin [10], martingale approximations have been an effective tool for studying the (conditional) central limit question and law of the iterated logarithm for stationary processes; see, for example, [3, 4, 20, 21], and their references for recent developments. The terminology here differs slightly from that of [20].

The sequence $X_k = g(W_k)$ is said to admit a *co-boundary* if there is a stationary sequence of martingale differences $d_k$ and another stationary process $Z_k$ for which

$$X_k = d_k + Z_k - Z_{k-1},$$

for all $k$, in which case $S_n = \tilde{M}_n + \tilde{R}_n$ with $\tilde{M}_n = d_1 + \cdots + d_n$ and $\tilde{R}_n = Z_n - Z_0$. Here $\tilde{M}_n$ is a martingale and $\tilde{R}_n$ is stochastically bounded, but does not necessarily satisfy $\|\tilde{R}_n\| = o(\sqrt{n})$. Conversely, a martingale approximation does not require $R_n$ to be stochastically bounded. The relation between co-boundaries and martingale approximations is further clarified by the examples of [9].

Letting $Q^*$ denote the adjoint of the restriction of $Q$ to $L^2(\pi)$, so that $\langle Qf, g \rangle = \langle f, Q^*g \rangle$ for $f, g \in L^2(\pi)$, $Q$ is said to be a co-isometry if $QQ^* = I$, in which case $Q^*$ is an isometry. Importantly, this condition is satisfied by shift processes. In Section 3, a convenient orthonormal basis for $L_0^2(\pi)$ is identified when $Q$ is a co-isometry, and a simple necessary and sufficient condition for the existence of a martingale approximation is given in terms of the coefficients in the expansion of $g$ with respect to this basis.

Returning to the main question, define

$$V_n g = \sum_{k=0}^{n-1} Q^k g,$$

so that $E(S_n | \mathcal{F}_1) = V_n g(W_1)$. If (1) holds, then $\|V_n g\|^2 = E[E(S_n | \mathcal{F}_1)^2] \leq 2E(M_1^2) + 2E(R_n^2) = o(n)$, and $\lim_{n \to \infty} E(S_n^2)/n = E(M_1^2)$. So, obvious necessary conditions for (1) are that

$$(2) \qquad \qquad \|V_n g\| = o(\sqrt{n})$$

and

$$(3) \qquad \qquad \|g\|_+^2 := \limsup_{n \to \infty} \frac{1}{n} E[S_n(g)^2] < \infty.$$

Let $\mathcal{L}$ denote the set of $g \in L_0^2(\pi)$ for which $\|g\|_+ < \infty$. Then $\mathcal{L}$ is a linear space, and $\|\cdot\|_+$ is a pseudo norm on $\mathcal{L}$, called the *plus norm* below. Moreover, $Q$ maps $\mathcal{L}$ into itself, since

$$(4) \quad S_n(g) = S_n(Qg) + \sum_{k=1}^{n} [g(W_k) - Qg(W_{k-1})] + Qg(W_0) - Qg(W_n);$$

and, therefore, $\|Qg\|_+ \leq \|g\|_+ + \sqrt{E\{[g(W_1) - Qg(W_0)]^2\}}$. In Section 4 it is shown that $g$ admits a martingale approximation iff (2) holds and

$$\lim_{m \to \infty} \frac{1}{m} \sum_{k=1}^{m} \|Q^k g\|_+^2 = 0.$$

These results are used in Section 5 to study the relationship between martingale approximations and solutions to the fractional Poisson equation, $g = \sqrt{(I - Q)}h$. The relation between martingale approximations and the conditional central limit theorem is explored in Section 6 with special attention to superpositions of linear processes. Section 2 contains some preliminaries.

**2. Preliminaries.** In this section, upon exhibiting some preliminary facts, we establish a useful criterion for martingale approximations; and in particular, we show martingale approximations are unique. Let

$$\bar{V}_n = \frac{V_1 + \cdots + V_n}{n} = \sum_{k=0}^{n-1} \left(1 - \frac{k}{n}\right) Q^k.$$

Then

$$(5) \qquad E[S_n(g)^2] = 2n\langle g, \bar{V}_n g\rangle - n\|g\|^2,$$

from [2], page 219, and

$$(6) \qquad \bar{V}_n = \bar{V}_n Q^i + V_i - \frac{1}{n}QV_nV_i$$

for all $n \geq 1$, $i \geq 1$ by simple algebra and induction. Next, let $\pi_1$ denote the joint distribution of $W_0$ and $W_1$, define

$$(7) \qquad H_n(w_0, w_1) = V_n g(w_1) - QV_n g(w_0)$$

and $\bar{H}_n(w_0, w_1) = \bar{V}_n g(w_1) - Q\bar{V}_n g(w_0)$ for $w_0, w_1 \in \mathcal{W}$. Then $H_n$ and $\bar{H}_n$ are in $L^2(\pi_1)$.

LEMMA 1. *If (2) holds, then $S_k = M_{nk} + R_{nk}$ where $M_{nk} = \bar{H}_n(W_0, W_1) + \cdots + \bar{H}_n(W_{k-1}, W_k)$ and $\max_{k \leq n} \|R_{nk}\| = o(\sqrt{n})$.*

PROOF. The lemma is almost a special case of Theorem 1 of [20]. Using (6) with $i = 1$,

$$R_{nk} = S_k - M_{nk} = Q\bar{V}_n g(W_0) - Q\bar{V}_n g(W_k) + \frac{1}{n}S_k(QV_n g),$$

from which it follows that $\max_{k \leq n} \|R_{nk}\| \leq 3\max_{k \leq n}\|V_k g\|$, which is $o(\sqrt{n})$ by (2). □

Of course, $M_{nk}$ is a martingale in $k$ for each $n$. The following proposition is closely related to Theorem 1 of [17].

PROPOSITION 1. *$g \in L_0^2(\pi)$ admits a martingale approximation iff (2) holds and $\bar{H}_n$ converges to a limit $H$ in $L^2(\pi_1)$, in which case*

$$(8) \qquad M_n = M_n(g) := \sum_{k=1}^n H(W_{k-1}, W_k).$$

*Consequently, martingale approximations are unique.*

PROOF. Suppose first that $g$ admits a martingale approximation, $S_n = M_n + R_n$. Then (2) holds and $S_n = M_{nn} + R_{nn}$, where $\|R_{nn}\| = o(\sqrt{n})$, by Lemma 1. So,

$$nE\{[\bar{H}_n(W_0, W_1) - M_1]^2\} = E[(M_{nn} - M_n)^2] = E[(R_{nn} - R_n)^2] = o(n),$$

implying the convergence of $\bar{H}_n(W_0, W_1)$ in $L^2(P)$; and this is equivalent to the convergence of $\bar{H}_n$ in $L^2(\pi_1)$.

Conversely, if (2) holds and $\bar{H}_n$ converges to a limit $H$, say; we can let $M_n = H(W_0, W_1) + \cdots + H(W_{n-1}, W_n)$ and $R_n = S_n - M_n$. Then (1) holds, $R_n = M_{nn} - M_n + R_{nn}$, and $\|R_n\| \le \sqrt{n}\|\bar{H}_n - H\| + \|R_{nn}\| = o(\sqrt{n})$, establishing both the sufficiency and (8). That martingale approximations are unique is then clear. □

COROLLARY 1. *If $g$ admits a martingale approximation, then so does $Q^k g$, and $M_1(Q^k g) = H(W_0, W_1) - H_k(W_0, W_1)$ with $H_k$ as defined in (7).*

PROOF. For $k = 1$, this follows directly from (4); and for $k = 2, 3, \ldots$, it follows by induction. □

As a second corollary, we may obtain necessary and sufficient conditions for a linear process. Let $\ldots, \xi_{-1}, \xi_0, \xi_1, \ldots$ be i.i.d. random variables with mean 0 and unit variance; let $a_0, a_1, a_2, \ldots$ be a square summable sequence; and consider a causal linear process

$$(9) \qquad X_j = \sum_{i=0}^{\infty} a_i \xi_{j-i} = \sum_{i \le j} a_{j-i} \xi_i.$$

Such a process is of the form $X_k = g(W_k)$, where $W_k = (\ldots, \xi_{k-1}, \xi_k)$. Letting $b_{-1} = 0$, $b_n = a_0 + \cdots + a_n$ for $n \ge 0$, and using (9),

$$S_n = \sum_{i \le 1}(b_{n-i} - b_{-i})\xi_i + \sum_{i=0}^{n-2} b_i \xi_{n-i},$$

where the first term on the right-hand side is $E(S_n|W_1)$. It follows that

$$\|V_n g\|^2 = \|E(S_n|W_1)\|^2 = \sum_{i=-1}^{\infty}(b_{i+n} - b_i)^2,$$

also, $V_n g(W_1) - QV_n g(W_0) = b_n \xi_1$, and $\bar{H}_n(W_0, W_1) = \bar{b}_n \xi_1$ with $\bar{b}_n = (b_1 + \cdots + b_n)/n$. Thus, for a linear process, (2) specializes to

$$(10) \qquad \lim_{n \to \infty} \frac{1}{n} \sum_{i=-1}^{\infty}(b_{i+n} - b_i)^2 = 0.$$

COROLLARY 2. *For the linear process defined in (9), the following are equivalent:*

(a) *There is a martingale approximation.*
(b) *Equation (10) holds and $\bar{b}_n$ converges.*
(c) *Equation (10) holds and $\bar{b}_n^2$ converges.*

Proof. In this case $\|\bar{H}_n - \bar{H}_m\|^2 = (\bar{b}_n - \bar{b}_m)^2$. Hence, (a) and (b) are equivalent by Proposition 1. It is clear that (b) implies (c) and it remains only to show that (c) implies (b). If $\bar{b}_n^2$ converges, but $\bar{b}_n$ does not, then $\bar{b}_n$ would have to oscillate between two values, there would be a positive $\varepsilon$ for which $|\bar{b}_{n+1} - \bar{b}_n| \geq \varepsilon$ infinitely often; but this is impossible, since $\bar{b}_{n+1} - \bar{b}_n = (b_{n+1} - \bar{b}_n)/(n+1)$ and $b_n = O(\sqrt{n})$, as $a_0, a_1, \ldots$ are square summable. □

In the next section, we show how to extend this example from linear functions of shift processes to measurable ones with mean 0 and finite variance.

**3. Co-isometries.** We suppose throughout this section that the chain has a trivial left tail field and that $Q$ is a co-isometry; that is,

(11)
$$\lim_{n \to \infty} \|Q^n f\| = 0 \quad \text{and} \quad QQ^* = I$$

for all $f \in L_0^2(\pi)$. We also suppose $L_0^2(\pi)$ is separable. These conditions are satisfied, for example, by (one-sided) shift processes.

With a view toward later examples, we work with $\mathscr{L}_0^2(\pi)$, the space of complex-valued, square integrable functions with mean 0 under $\pi$. Then (11) is still valid for this space if we extend the definition of $Q$ to the imaginary part.

Let $\mathcal{H}$ denote a closed linear subspace of $\mathscr{L}_0^2(\pi)$ that is invariant under both $Q$ and $Q^*$; restrict $Q$ and $Q^*$ to $\mathcal{H}$; and let $\mathcal{K} = Q^*\mathcal{H}$. Then $Q^*$ is an isometry from $\mathcal{H}$ onto $\mathcal{K}$, since $\langle Q^*f, Q^*g \rangle = \langle f, QQ^*g \rangle = \langle f, g \rangle$ for $f, g \in \mathcal{H}$. This is the origin of the term "co-isometry." Moreover,

(12)
$$Q^*h(W_1) = h(W_0) \qquad \text{w.p.1}$$

for any $h \in L_0^2(\pi)$, since $E[Q^*h(W_1)h(W_0)] = \langle QQ^*h, h \rangle = \|h\|^2$ by conditioning on $W_0$, and therefore, $E\{[Q^*h(W_1) - h(W_0)]^2\} = \|Q^*h\|^2 - 2\langle QQ^*h, h \rangle + \|h\|^2 = 0$. It can be easily checked (12) also holds for $h \in \mathscr{L}_0^2(\pi)$.

Lemma 2. $\mathcal{K}$ is a closed, proper linear subspace of $\mathcal{H}$; and $\bigcap_{j=0}^{\infty} Q^{*j}\mathcal{H} = \{0\}$.

Proof. That $\mathcal{K}$ is closed is clear, since $Q^*$ is an isometry; and that $\mathcal{K}$ is proper follows from $\bigcap_{j=0}^{\infty} Q^{*j}\mathcal{H} = \{0\}$. So, it suffices to establish the latter. If $f \in \bigcap_{j=0}^{\infty} Q^{*j}\mathcal{H}$, then there are $h_0, h_1, \ldots \in \mathcal{H}$ for which $f = Q^{*j}h_j$ with each $j$. In this case, $\|h_j\| = \|f\|$, since $Q^*$ is an isometry, $h_j = Q^j Q^{*j}h_j = Q^j f$, and $\lim_{j \to \infty} \|Q^j f\| = 0$. So, $\|f\| = 0$, establishing the lemma. □

Next, let $\mathcal{K}^{\perp} = \{f \in \mathcal{H} : \langle f, h \rangle = 0 \text{ for all } h \in \mathcal{K}\}$. Then $\mathcal{K}^{\perp} = \{g \in \mathcal{H} : Qg = 0\}$, since $\langle Q^*f, g \rangle = \langle f, Qg \rangle = 0$ for all $f \in \mathcal{H}$ iff $Qg = 0$; and $Q^*Q$ is the projection operator onto $\mathcal{K}$, since $(Q^*Q)^2 = Q^*Q$ and $Q(I - Q^*Q) = 0$. Let $\mathcal{E}_0 = \{e_j : j \in J\}$ be an orthonormal basis for $\mathcal{K}^\perp$, let $\mathcal{E}_i = Q^{*i}\mathcal{E}_0$ and $\mathcal{E} = \bigcup_{i \geq 0} \mathcal{E}_i$.

LEMMA 3. $\mathcal{E}$ is an orthonormal basis for $\mathcal{H}$.

PROOF. $\mathcal{E}_i$ consists of orthonormal elements for each $i \geq 0$, since $Q^*$ is an isometry; for any $f \in \mathcal{E}_i$ and $f' \in \mathcal{E}_{i'}$, where $i < i'$, there are $e, e' \in \mathcal{E}_0$ for which $f = Q^{*i}e$ and $f' = Q^{*i'}e'$, in which case $\langle f, f' \rangle = \langle Q^{*i}e, Q^{*i'}e' \rangle = \langle Q^{i'-i}e, e' \rangle = 0$, since $Qe = 0$. Finally, if $f \perp \mathcal{E}_0$, then $f \in \mathcal{K}$ and $f = Q^*h_1$ for some $h_1 \in \mathcal{H}$. If also, $f \perp Q^*\mathcal{E}_0$, then $Qf \perp \mathcal{E}_0$, $Qf = Q^*h_2$ for some $h_2 \in \mathcal{H}$, and $f = Q^*Qf = Q^{*2}h_2$. Continuing, we find that if $f \perp \mathcal{E}$, then $f \in Q^{*j}\mathcal{H}$ for all $j$, and completeness follows from Lemma 2. $\square$

Now write $e_{i,j} = Q^{*i}e_j$, so that $\mathcal{E}_i = \{e_{i,j} : j \in J\}$, and let $\mathcal{H}_j = \mathrm{span}(e_{i,j} : i \geq 0)$, the closed linear span of $\{e_{i,j} : i \geq 0\}$. Then $Q\mathcal{H}_j = \mathcal{H}_j$ for each $j$, and $\mathcal{H} = \bigoplus_{j \in J} \mathcal{H}_j$. In the language of [5, 11], the $\mathcal{H}_j$, $j \in J$, are an orthogonal invariant splitting of $\mathcal{H}$. Then, any $g \in \mathcal{H}$ may be written as $g = \sum_{j \in J} \sum_{i=0}^{\infty} c_{i,j} e_{i,j}$, where $c_{i,j}$ are square summable. Let $b_{n,j} = c_{0,j} + \cdots + c_{n-1,j}$, $\bar{b}_{n,j} = (b_{1,j} + \cdots + b_{n,j})/n$ and regard $\mathbf{b}_n = (b_{n,j} : j \in J)$ and $\bar{\mathbf{b}}_n = (\bar{b}_{n,j} : j \in J)$ as elements of $\ell^2(J)$.

THEOREM 1. $g \in L_0^2(\pi)$ admits a martingale approximation iff $\bar{\mathbf{b}}_n$ converges in $\ell^2(J)$, and

$$(13) \qquad \lim_{n \to \infty} \frac{1}{n} \sum_{i=0}^{\infty} \|\mathbf{b}_{i+n} - \mathbf{b}_i\|^2 = 0.$$

PROOF. We take $\mathcal{H} = \mathscr{L}_0^2(\pi)$. Since $Qe_{i,j} = QQ^{*i}e_j = 0$ if $i = 0$ and $e_{i-1,j}$ if $i \geq 1$, $Qg = \sum_{j \in J} \sum_{i=1}^{\infty} c_{i,j} e_{i-1,j}$,

$$Q^k g = \sum_{i=k}^{\infty} \sum_{j \in J} c_{i,j} e_{i-k,j} = \sum_{i=0}^{\infty} \sum_{j \in J} c_{i+k,j} e_{i,j},$$

$$V_n g = \sum_{i=0}^{\infty} \sum_{j \in J} (b_{i+n,j} - b_{i,j}) e_{i,j}$$

and

$$\|V_n g\|^2 = \sum_{i=0}^{\infty} \sum_{j \in J} |b_{i+n,j} - b_{i,j}|^2 = \sum_{i=0}^{\infty} \|\mathbf{b}_{i+n} - \mathbf{b}_i\|^2.$$

So (13) is just (2), specialized to the present context.

Next $Q^*Qg = \sum_{j \in J} \sum_{i=1}^{\infty} c_{i,j} e_{i,j}$, so that from (12),

$$g(W_1) - Qg(W_0) = [g - Q^*Qg](W_1) = \sum_{j \in J} c_{0,j} e_{0,j}(W_1),$$

$$\bar{H}_n(W_0, W_1) = \sum_{k=0}^{n-1} \left(1 - \frac{k}{n}\right)[Q^k g(W_1) - Q^{k+1} g(W_0)] = \sum_{j \in J} \bar{b}_{n,j} e_{0,j}(W_1)$$

and

$$\|\bar{H}_n - \bar{H}_m\| = \|\bar{\mathbf{b}}_n - \bar{\mathbf{b}}_m\|.$$

The theorem now follows directly from Proposition 1.  □

EXAMPLE 1 (Bernoulli shifts).   The one-sided Bernoulli shift process is defined by

$$W_k = \sum_{j=0}^{\infty} (\tfrac{1}{2})^{j+1} \xi_{k-j},$$

where $\ldots, \xi_{-1}, \xi_0, \xi_1, \ldots$ are i.i.d. random variables taking the values 0 and 1 with probability 1/2 each. The state space $\mathcal{W}$ is the unit interval, the marginal distribution $\pi$ is the uniform distribution, $Qg(w) = \frac{1}{2}[g(\frac{1}{2}w) + g(\frac{1}{2}w + \frac{1}{2})]$, and $Q^*g(w) = g(2w)$ for a.e. $w \in \mathcal{W}$ and $g \in L^1(\pi)$ with the convention that $g$ is continued periodically. For this example, any $g \in \mathscr{L}_0^2(\pi)$ has a Fourier expansion

(14) $$g = \sum_{r \neq 0} c_r e_r,$$

where $e_r(w) = e^{2\pi \imath r w}$ and $c_r$, $r \in \mathbb{Z}$, are square summable. Then $Qe_r = 0$ or $e_{(1/2)r}$ accordingly as $r$ is odd or even, and $Q^*e_r = e_{2r}$ for all $r$. With $\mathcal{H} = \mathscr{L}_0^2(\pi)$, it follows that $\mathcal{K}$, respectively $\mathcal{K}^\perp$, consists of all functions $g$ for which $c_r = 0$ for odd, respectively even, $r$. Thus, $\mathcal{E}_0 = \text{span}(e_r : r \in \text{Odd})$, and $\mathcal{E}_i = \text{span}(e_{r2^i} : r \in \text{Odd})$, and there is an invariant splitting with $e_{i,j} = e_{j2^i}$. Necessary and sufficient conditions for the existence of a martingale approximation can be read from Theorem 1. See [19] for more on the Fourier analysis of Bernoulli shifts.

EXAMPLE 2 (Lebesgue shifts).   By a (one-sided) Lebesgue shift, we mean the Markov chain $W_k = (\ldots, U_{k-1}, U_k)$ where $\ldots, U_{-1}, U_0, U_1, \ldots$ are independent uniformly distributed random variables over $[0, 1)$, in which case $\mathcal{W} = [0, 1)^{\mathbb{N}}$ and $\pi = \lambda^{\mathbb{N}}$, where $\lambda$ is the uniform distribution. Lebesgue shifts are similar to Bernoulli shifts. Let $\Gamma$ denote the set of sequences $j = (j_0, j_1, \ldots) \in \mathbb{Z}^{\mathbb{N}}$ for which $j_i = 0$ for all but finite number of $i$. Then, letting $j \cdot w = j_0 w_0 + j_1 w_{-1} + \cdots$ and $e_j(w) = e^{2\pi \imath j \cdot w}$ for $w = (\ldots, w_{-1}, w_0) \in [0, 1)^{\mathbb{N}}$ and $j \in \Gamma$, any $g \in \mathscr{L}_0^2(\pi)$ has a Fourier expansion,

$$g(w) = \sum_{j \in \Gamma} c_j e_j$$

where $c_j$ are square summable. Next, let $J = \{j \in \Gamma : j_0 \neq 0\}$. Then, since

$$Qe_j(w) = \left[ \int_0^1 e^{2\pi \imath j_0 u} \, du \right] \prod_{i=1}^{\infty} \exp(2\pi \imath j_i w_{-i+1}),$$

$\mathcal{E}_0 = \{e_k : k \in J\}$ is an orthonormal basis for $\mathcal{K}^{\perp}$ [with $\mathcal{H} = \mathscr{L}_0^2(\pi)$ and $\mathcal{K} = Q^*\mathcal{H}$]. Define $\psi \colon \Gamma \to \Gamma$ by $\psi(j) = (0, j_0, j_1, \ldots)$, then it is not difficult to check $Q^* e_k = e_{\psi(k)}$, $Q^{*i} e_k = e_{\psi^i(k)}$, where $\psi^i$ is the composition of $\psi$ with itself $i$ times. Necessary and sufficient conditions can be read from Theorem 1.

EXAMPLE 3 (Superlinear processes). Let $\xi_{i,j}$, $i \in \mathbb{Z}$, $j \in \mathbb{N}$, be independent random variables, all having mean 0 and bounded variances, for which $\ldots, \xi_{-1,j}, \xi_{0,j}, \xi_{1,j}, \ldots$ are identically distributed for each $j$, and let $c_{ij}$, $i \in \mathbb{Z}$, $j \in \mathbb{N}$, be a square summable array. Then

$$(15) \qquad X_k = \sum_{j=0}^{\infty} \sum_{i=0}^{\infty} c_{i,j} \xi_{k-i,j}$$

converges w.p.1 and in mean square for each $k$ and defines a stationary process. Letting $\boldsymbol{\xi}_i = (\xi_{i,0}, \xi_{i,1}, \ldots)$, $X_k$ is of the form $X_k = g(W_k)$, where $W_k = (\ldots, \boldsymbol{\xi}_{k-1}, \boldsymbol{\xi}_k)$ is a shift process. Next, letting $\mathcal{H} = \operatorname{span}(\xi_{i,j} : i \leq 0, j \geq 0)$, one finds easily that there is an invariant splitting with $e_{i,j} = \xi_{-i,j}$ for $i, j \geq 0$. Necessary and sufficient conditions for the existence of a martingale approximation can again be read from Theorem 1.

**4. The plus norm.** To study the plus norm, we first recall the definition $\|g\|_+^2 = \limsup_{n \to \infty} E[S_n(g)^2]/n$. The following example serves as a simple illustration.

EXAMPLE 4. If $Q$ is a co-isometry and the chain has a trivial left tail field, we may write $g = \sum_{j \in J} \sum_{i=0}^{\infty} c_{i,j} e_{i,j}$, as in Section 3, and $\bar{H}_n(W_0, W_1) = \sum_{j \in J} \bar{b}_{n,j} e_{0,j}(W_1)$, as in the proof of Theorem 1. So, if (2) holds, $E(S_n^2) = nE[\bar{H}_n^2(W_0, W_1)] + o(n) = n\|\bar{\mathbf{b}}_n\|^2 + o(n)$, and $\|g\|_+^2 = \limsup_{n \to \infty} \|\bar{\mathbf{b}}_n\|^2$.

The main result of this section is that $g$ admits a martingale approximation iff $\|V_n g\| = o(\sqrt{n})$ and $\sum_{k=1}^{m} \|Q^k g\|_+^2 = o(m)$. The following two lemmas are needed; their proofs are given after the proof of Theorem 2.

LEMMA 4.   *If $g \in L_0^2(\pi)$ and (2) holds, then*

$$\lim_{n \to \infty} \left[ \| \bar{H}_n - \bar{H}_m \|^2 - \frac{2}{m} \langle \bar{V}_n g, QV_m g \rangle \right] = -\frac{2}{m} \langle \bar{V}_m g, QV_m g \rangle - \left\| \frac{QV_m g}{m} \right\|^2.$$

LEMMA 5.   *If $g \in L_0^2(\pi)$ and $\|g\|_+ < \infty$, then*

$$\limsup_{n \to \infty} \langle \bar{V}_n g, QV_m g \rangle \leq \tfrac{1}{2} \sum_{k=1}^m \|Q^k g\|_+^2 + \tfrac{1}{2} \|QV_m g\|^2 + \langle g, V_m Q g \rangle;$$

*and if $g$ admits a martingale approximation, then the limit exists and there is equality.*

THEOREM 2.   *$g$ admits a martingale approximation iff (2) holds and*

$$(16) \qquad\qquad \lim_{m \to \infty} \frac{1}{m} \sum_{k=1}^m \|Q^k g\|_+^2 = 0.$$

PROOF.   Suppose first that $g$ admits a martingale approximation. Then $\|V_n g\| = o(\sqrt{n})$ and $\lim_{m \to \infty}[\lim_{n \to \infty} \|\bar{H}_n - \bar{H}_m\|^2] = 0$ by Proposition 1. Next, by Lemmas 4 and 5,

$$\begin{aligned}
\lim_{n \to \infty} \| \bar{H}_n - \bar{H}_m \|^2 &= \lim_{n \to \infty} \frac{2}{m} \langle \bar{V}_n g, QV_m g \rangle - \left[ \frac{2}{m} \langle \bar{V}_m g, QV_m g \rangle + \left\| \frac{QV_m g}{m} \right\|^2 \right] \\
&= \frac{1}{m} \sum_{k=1}^m \|Q^k g\|_+^2 + \frac{1}{m} \|QV_m g\|^2 + \frac{2}{m} \langle g, QV_m g \rangle \\
&\quad - \frac{2}{m} \langle \bar{V}_m g, QV_m g \rangle - \left\| \frac{QV_m g}{m} \right\|^2.
\end{aligned}$$

Since $\|V_m g\| = o(\sqrt{m})$, the last four terms on the right approach 0 as $m \to \infty$, and, therefore, so does the first. This establishes the necessity of (16).

Next suppose that (2) and (16) hold; then $\lim_{m \to \infty}[\limsup_{n \to \infty} \|\bar{H}_n - \bar{H}_m\|^2] = 0$, by Lemmas 4 and 5. It follows easily that $\sup_{n \geq 1} \|\bar{H}_n\| < \infty$, which implies $\bar{H}_1, \bar{H}_2, \ldots$ is weakly compact in $L^2(\pi_1)$. Let $\bar{H}^*$ denote any weak limit point of $\bar{H}_1, \bar{H}_2, \ldots$. Then $\|H^* - \bar{H}_m\| \leq \limsup_{n \to \infty} \|\bar{H}_n - \bar{H}_m\|$ for each $m$ (cf. [8], page 68). Thus, $\lim_{m \to \infty} \|\bar{H}_m - H^*\| = 0$ from which the converse follows from Proposition 1.   □

PROOF OF LEMMA 4.   To begin, write

$$\begin{aligned}
\| \bar{H}_n - \bar{H}_m \|^2 &= \|(\bar{V}_n - \bar{V}_m)g\|^2 - \|Q(\bar{V}_n - \bar{V}_m)g\|^2 \\
&= \langle (I+Q)(\bar{V}_n - \bar{V}_m)g, (I-Q)(\bar{V}_n - \bar{V}_m)g \rangle \\
&= 2 \left\langle (\bar{V}_n - \bar{V}_m)g, \left( \frac{QV_m g}{m} - \frac{QV_n g}{n} \right) \right\rangle - \left\| \frac{QV_m g}{m} - \frac{QV_n g}{n} \right\|^2;
\end{aligned}$$

and when the first term in the last line is expanded, it becomes

$$\frac{2}{m}\langle \bar{V}_n g, QV_m g\rangle - \frac{2}{n}\langle \bar{V}_n g, QV_n g\rangle - \frac{2}{m}\langle \bar{V}_m g, QV_m g\rangle + \frac{2}{n}\langle \bar{V}_m g, QV_n g\rangle.$$

The lemma now follows directly from (2) and the mean ergodic theorem, which implies that all those terms multiplied by $1/n$ approach 0 as $n \to \infty$. □

PROOF OF LEMMA 5. Writing

$$\langle \bar{V}_n g, QV_m g\rangle = \sum_{k=1}^{m}\langle \bar{V}_n g, Q^k g\rangle,$$

and using (6), then

$$\langle \bar{V}_n g, QV_m g\rangle = \sum_{k=1}^{m}\left[\langle \bar{V}_n Q^k g, Q^k g\rangle + \langle V_k g, Q^k g\rangle - \frac{1}{n}\langle QV_n V_k g, Q^k g\rangle\right].$$

Here

$$\sum_{k=1}^{m}\langle V_k g, Q^k g\rangle = \sum_{k=1}^{m}\sum_{j=1}^{k}\langle Q^{j-1} g, Q^k g\rangle$$

$$= \frac{1}{2}\sum_{k=1}^{m}\sum_{j=1}^{m}\langle Q^j g, Q^k g\rangle - \frac{1}{2}\sum_{j=1}^{m}\|Q^j g\|^2 + \langle g, V_m Qg\rangle$$

$$= \frac{1}{2}\|V_m Qg\|^2 - \frac{1}{2}\sum_{j=1}^{m}\|Q^j g\|^2 + \langle g, V_m Qg\rangle.$$

Combining terms together,

$$\langle \bar{V}_n g, QV_m g\rangle = \frac{1}{2}\sum_{k=1}^{m}[2\langle \bar{V}_n Q^k g, Q^k g\rangle - \|Q^k g\|^2]$$

$$+ \frac{1}{2}\|QV_m g\|^2 + \langle g, V_m Qg\rangle - \sum_{k=1}^{m}\frac{1}{n}\langle QV_n V_k g, Q^k g\rangle.$$

The first assertion follows directly from (2) and (5). So does the second; for if $g$ admits a martingale approximation, then the limit exists in the definition of $\|Q^k g\|_+$. □

**5. The fractional Poisson equation.** It is possible to attach a meaning to the symbol $\sqrt{I-Q}$ by replacing $t$ with $Q$ in the series expansion of $\sqrt{1-t}$. The definition may be written

$$\sqrt{I-Q} = I - \sum_{k=1}^{\infty}\beta_k Q^k,$$

where $\beta_k = (-1)^{k-1} \binom{1/2}{k}$ and the series converges in the operator norm, since $\beta_k \sim 1/(2\sqrt{\pi}k^{3/2})$ as $k \to \infty$. A function $h \in L_0^2(\pi)$ is said to solve the *fractional Poisson equation* (for $g$) if $g = \sqrt{(I-Q)}h$. The relation between the existence of a solution to the fractional Poisson equation and the existence of a martingale approximation is considered in this section for co-isometries and normal operators ($QQ^* = Q^*Q$).

LEMMA 6.  *If $g \in \sqrt{(I-Q)}L_0^2(\pi)$, then $\|V_n g\| = o(\sqrt{n})$; and if $g = \sqrt{(I-Q)}h = \sqrt{(I-Q^*)}h^*$, then $\|g\|_+^2 = \langle (I+Q)h, h^* \rangle$.*

PROOF.  Observe that $(I-Q^k)V_n = (I-Q^n)V_k$. So, if $g = \sqrt{(I-Q)}h$, then $V_n g = \sum_{k=0}^\infty \beta_k (I - Q^{k \vee n}) V_{k \wedge n} h$, where $\wedge$ ($\vee$) denotes minimum (maximum). Using the mean ergodic theorem, $\|V_n h\| = o(n)$, then

$$\|V_n g\| \le 2 \sum_{k=0}^\infty \beta_k \|V_{k \wedge n} h\| = 2 \sum_{k=0}^\infty \beta_k o(k \wedge n) = o(\sqrt{n}),$$

establishing the first assertion. If, in addition, $g = \sqrt{(I-Q^*)}h^*$, then $\|g\|^2 = \langle (I-Q)h, h^* \rangle$, and

$$\langle \bar{V}_n g, g \rangle = \langle (I-Q)\bar{V}_n h, h^* \rangle = \langle h, h^* \rangle - \frac{1}{n} \langle Q V_n h, h^* \rangle \to \langle h, h^* \rangle,$$

using the mean ergodic theorem again in the final step. Thus, in view of (5), $\|g\|_+^2 = \lim_{n \to \infty} [2 \langle \bar{V}_n g, g \rangle - \|g\|^2] = \langle (I+Q)h, h^* \rangle$; similar calculations also appear in [7].  □

*Normal operators.*  As an interesting generalization of [13] in the reversible case, it is known [1, 6, 12] that if $Q$ is a normal operator and there is a solution to the fractional Poisson equation, then $g$ admits a martingale approximation. This result can be easily deduced from our Theorem 2. Recall that if $R$ is any bounded normal operator on a Hilbert space $\mathcal{H}$, then $\sqrt{I-R}$ and $\sqrt{I-R^*}$ have the same range (cf. [6], Lemma 2).

PROPOSITION 2.  *Suppose that $Q$ is normal; then any $g \in \sqrt{(I-Q)}L_0^2(\pi)$ admits a martingale approximation.*

PROOF.  If $g \in \sqrt{(I-Q)}L_0^2(\pi)$, then (2) follows from Lemma 6, and it suffices to establish (16). Since the ranges of $\sqrt{(I-Q)}$ and $\sqrt{(I-Q^*)}$ are the same, there are $h, h^* \in L_0^2(\pi)$ for which $g = \sqrt{(I-Q)}h = \sqrt{(I-Q^*)}h^*$. Then $Q^k g = \sqrt{(I-Q)}Q^k h = \sqrt{(I-Q^*)}Q^k h^*$, so that $\|Q^k g\|_+^2 = \langle (I+Q)Q^k h, Q^k h^* \rangle$. Thus, letting $R = Q^*Q$, $\|Q^k g\|_+^2 = \langle (I+Q)h, R^k h^* \rangle$, and it is necessary to show

$$\lim_{m \to \infty} \frac{1}{m} \sum_{k=1}^m \langle (I+Q)h, R^k h^* \rangle = 0. \tag{17}$$

To see this let $\mathcal{R}$ be the closure of $(I - R)L_0^2(\pi)$. Then $\mathcal{R}^\perp$ consists of all $f$ for which $Rf = f$, and $Q$, $Q^*$, and $R$ map both $\mathcal{R}$ and $\mathcal{R}^\perp$ into themselves. Write $h = h_1 + h_2$ with $h_1 \in \mathcal{R}$, $h_2 \in \mathcal{R}^\perp$, and let $g_i = \sqrt{(I - Q)}h_i$. Then $g_1 \in \mathcal{R}$ and $g_2 \in \mathcal{R}^\perp$, since $Q$ maps $\mathcal{R}$ and $\mathcal{R}^\perp$ into themselves. Next, write $h^* = h_1^* + h_2^*$ with $h_1^* \in \mathcal{R}$, $h_2^* \in \mathcal{R}^\perp$; then $g_i = \sqrt{(I - Q^*)}h_i^*$ by the uniqueness of direct sum decomposition of $g$. Returning to (17), we have

$$\langle (I + Q)h, R^k h^* \rangle = \langle (I + Q)h_1, R^k h_1^* \rangle + \langle (I + Q)h_2, h_2^* \rangle$$
$$= \langle (I + Q)h_1, R^k h_1^* \rangle + \|g_2\|_+^2$$

by orthogonality and Lemma 6. It will be first shown that $\|g_2\|_+ = 0$; to see it, note $Rg_2 = g_2$, then

$$\|V_n g_2\|^2 = 2 \sum_{j=0}^{n-1} \sum_{k=j}^{n-1} \langle Q^j g_2, Q^k g_2 \rangle - \sum_{j=0}^{n-1} \|Q^j g_2\|^2$$
$$= 2 \sum_{j=0}^{n-1} \langle g_2, V_{n-j} g_2 \rangle - n \|g_2\|^2$$
$$= n[2\langle g_2, \bar{V}_n g_2 \rangle - \|g_2\|^2],$$

thus, $\|g_2\|_+ = 0$ follows from (5) and Lemma 6. That (17) holds when $h_1^* \in (I - R)L_0^2(\pi)$ is clear by forming a telescoping sum, and the boundary case then follows by approximation. □

*Co-isometries.* The existence of a solution to the fractional Poisson equation does not imply the existence of a martingale approximation for co-isometries. Here is a simple example.

EXAMPLE 5. Let $\ldots, \xi_{-1}, \xi_0, \xi_1, \ldots$ be i.i.d. with mean 0 and unit variance; consider the shift process $W_k = (\ldots, \xi_{k-1}, \xi_k)$. For $j \geq 0$, let $a_j = 1/[\sqrt{(j+1)} \log(j+2)]$ and define $h$ by

$$h(W_0) = \sum_{j=0}^{\infty} a_j \xi_{-j},$$

so that $h(W_k)$ is a linear process. Then $g = \sqrt{(I - Q)}h$ admits a solution to the fractional Poisson equation, and

$$g(W_0) = \sum_{k=1}^{\infty} \beta_k (I - Q^k)h(W_0) = \sum_{j=0}^{\infty} c_j \xi_{-j}$$

with

$$c_j = \sum_{k=1}^{\infty} \beta_k (a_j - a_{j+k}),$$

after some straightforward calculation. Observe that $a_j - a_{j+k} \geq 0$ for all $j \geq 0$ and $k \geq 1$, and that $a_{j+k} \leq 3a_j/4$ for all $k \geq j+1$ and all $j \geq 0$. So,

$$c_j \geq \frac{1}{4}a_j \sum_{k=j+1}^{\infty} \beta_k \geq \left(\frac{a_j}{9\sqrt{j}}\right)$$

for all sufficiently large $j$. Therefore, $b_n = c_0 + \cdots + c_n \to \infty$, and also, its Cesàro average $\bar{b}_n \to \infty$ as $n \to \infty$. No martingale approximation can exist.

However, the existence of solutions to both the forward and backward fractional Poisson equations does imply the existence of a martingale approximation.

PROPOSITION 3.  *Suppose $Q$ is a co-isometry and the chain has a trivial left tail field, and if $g \in \sqrt{(I-Q)}L_0^2(\pi) \cap \sqrt{(I-Q^*)}L_0^2(\pi)$, then $g$ admits a martingale approximation.*

PROOF.  As in Section 3, we can take $\mathcal{H} = \mathscr{L}_0^2(\pi)$, and there is an orthogonal invariant splitting, $\mathcal{H} = \bigoplus_{j \in J} \mathcal{H}_j$. Let $g = \sqrt{(I-Q)}h$ for some $h \in \mathcal{H}$, $g = \sum_{j \in J} g_j$, and $h = \sum_{j \in J} h_j$ with $g_j, h_j \in \mathcal{H}_j$ for all $j$. Clearly $g = \sum_{j \in J} \sqrt{(I-Q)}h_j$ and, therefore, $g_j = \sqrt{(I-Q)}h_j$, by taking the projection on each $\mathcal{H}_j$. Similarly, $g = \sqrt{(I-Q^*)}h^*$, where $h^* = \sum_{j \in J} h_j^*$ with $h_j^* \in \mathcal{H}_j$, and $g_j = \sqrt{(I-Q^*)}h_j^*$ for each $j$. It then follows easily from Lemma 6 and Example 4 that $\lim_{n \to \infty} |\bar{b}_{n,j}|^2 = \|g_j\|_+^2 = \langle (I+Q)h_j, h_j^* \rangle$ exists for each $j$ and that $\lim_{n \to \infty} \|\bar{\mathbf{b}}_n\|^2 = \|g\|_+^2 = \langle (I+Q)h, h^* \rangle$ exist. It then follows from (the proof of) Corollary 2 that $b_j = \lim_{n \to \infty} \bar{b}_{n,j}$ exists for each $j$, so that $\bar{\mathbf{b}}_n$ converges weakly to $\mathbf{b} = (b_j : j \in J)$. So, to show convergence of $\bar{\mathbf{b}}_n$ in the norm of $\ell^2(J)$ and, therefore, the existence of a martingale approximation, it suffices to show that $\lim_{n \to \infty} \|\bar{\mathbf{b}}_n\|^2 = \|\bar{\mathbf{b}}\|^2$; and this follows easily from Lemma 6 which implies

$$\lim_{n \to \infty} \|\bar{\mathbf{b}}_n\|^2 = \langle (I+Q)h, h^* \rangle = \sum_{j \in J} \langle (I+Q)h_j, h_j^* \rangle = \sum_{j \in J} |b_j|^2 = \|\mathbf{b}\|^2. \quad \square$$

**6. The CCLT for superlinear processes.**  Let $F_n$ denote the conditional distribution function of $S_n/\sqrt{n}$ given $W_0$,

$$F_n(w; z) = P\left[\frac{S_n}{\sqrt{n}} \leq z \Big| W_0 = w\right].$$

We will say that the *conditional central limit theorem* (CCLT) holds (with a $\sqrt{n}$ normalization) iff

$$\lim_{n \to \infty} \frac{E[S_n^2]}{n} = \kappa^2 \in [0, \infty)$$

and

$$\lim_{n \to \infty} \int_{\mathcal{W}} d[\Phi_\kappa, F_n(w; \cdot)] \pi\{dw\} = 0,$$

where $\Phi_\kappa$ denotes the normal distribution function with mean 0 and standard deviation $\kappa$, and $d$ is the Lévy metric or any other bounded metric that metrizes weak convergence of distribution functions.

It is clear that the existence of a martingale approximation implies the CCLT; see, for example, [15]. It is also clear, for simple linear process as defined in (9), CCLT necessarily requires the existence of martingale approximation. However, in general, the converse is not true as shown in the example below. To proceed as in Example 3, let $F_j$ be the common distribution function of $\xi_{i,j}$, $i = \ldots, -1, 0, 1, \ldots$ and suppose that the $F_j$ have mean 0 and bounded variances. Recall the notation $b_{n,j} = c_{0,j} + \cdots + c_{n-1,j}$ and $\bar{b}_{n,j} = (b_{1,j} + \cdots + b_{n,j})/n$ and that $\mathbf{b}_n = (b_{n,1}, b_{n,2}, \ldots)$ and $\bar{\mathbf{b}}_n$ may be regarded as elements of $\ell^2(\mathbb{N})$.

EXAMPLE 6 (Superlinear process revisited). Consider a superlinear process, defined in (15), with $c_{i,j} = 0$ for all $j \geq 2$, $b_{n,0} = \cos(\sqrt{\log n})$, $b_{n,1} = \sin(\sqrt{\log n})$, and $c_{0,j} = c_{1,j} = 0$ for $j = 0, 1$. Then $c_{n,j} = b_{n,j} - b_{n-1,j} = O(1/(n\sqrt{\log n}))$ for $j = 0, 1$. So, the process is well defined. If $F_0$ and $F_1$ both have mean 0 and unit variance, then the CCLT holds, but martingale approximation does not exist. To see this, first observe that for any $\delta > 0$,

$$\sum_{k=0}^{\infty} (b_{k+n,0} - b_{k,0})^2 \leq \left( \sum_{k \leq n\delta} + \sum_{k > n\delta} \right) [\cos(\sqrt{\log(n+k)}) - \cos(\sqrt{\log k})]^2$$

$$\leq 4n\delta + \sum_{k > n\delta} \left( \frac{n}{2k\sqrt{\log k}} \right)^2,$$

so that $\sum_{k=0}^{\infty} (b_{k+n,0} - b_{k,0})^2 = o(n)$, and similarly, $\sum_{k=0}^{\infty} (b_{k+n,1} - b_{k,1})^2 = o(n)$. So, $\|V_n g\|^2 = o(n)$. Next, for any $\varepsilon > 0$,

$$\bar{b}_{n,0} - b_{n,0} = \frac{1}{n} \sum_{k=1}^{n} [\cos(\sqrt{\log k}) - \cos(\sqrt{\log n})]$$

$$\leq \frac{1}{n} \sum_{k \leq n\varepsilon} [\cos(\sqrt{\log k}) - \cos(\sqrt{\log n})] + \frac{1}{n} \sum_{n\varepsilon < k \leq n} [\sqrt{\log n} - \sqrt{\log k}]$$

$$\leq 2\varepsilon + \frac{1}{n} \sum_{n\varepsilon < k \leq n} \frac{n-k}{2k\sqrt{\log k}} \leq 2\varepsilon + \frac{1}{n}(n - n\varepsilon) \frac{1}{2\varepsilon\sqrt{\log(n\varepsilon)}}$$

for all large $n$. It follows that $\bar{b}_{n,0} - b_{n,0} = o(1)$. Similarly $\bar{b}_{n,1} - b_{n,1} = o(1)$, and therefore, $\bar{b}_{n,0}^2 + \bar{b}_{n,1}^2 \to 1$. So, applying Theorem 2 of [20], CCLT holds; but martingale approximation does not exist since $\bar{b}_{n,j}$ does not converge for $j = 0, 1$.

Next, we investigate some partial converses for superlinear processes.

THEOREM 3. *If the CCLT holds for all choices $F_1, F_2, \ldots$ with means $0$ and unit variances, then $\bar{\mathbf{b}}_n$ is pre-compact in $\ell^2(\mathbb{N})$; and if the CCLT holds for all $F_1, F_2, \ldots$ with means $0$ and bounded variances, then $\bar{\mathbf{b}}_n$ converges in $\ell^2(\mathbb{N})$.*

PROOF. If the CCLT holds, then (2) holds by Corollary 1 of [15]. So, by Lemma 1, $S_n = M_{nn} + R_{nn}$, where $\|R_{nn}\| = o(\sqrt{n})$ and

$$M_{nn} = \sum_{j=1}^{\infty} \bar{b}_{n,j} \zeta_{n,j},$$

where $\zeta_{n,j} = \xi_{1,j} + \cdots + \xi_{n,j}$. So, if the CCLT holds for any choice of $F_1, F_2, \ldots$ with means $0$ and unit variances, then $\lim_{n \to \infty} \|\bar{\mathbf{b}}_n\|^2 = \kappa^2$. In particular, $\bar{\mathbf{b}}_n$, $n \geq 1$, are bounded and, therefore, weakly pre-compact. To show pre-compactness, it thus suffices to show that any weak limit point is a strong limit point. Let $\mathbf{b} \in \ell^2(\mathbb{N})$ be an arbitrary weak limit point and let $\mathbb{N}_0$ be a subsequence for which $\lim_{n \in \mathbb{N}_0} \bar{\mathbf{b}}_n = \mathbf{b}$. Then $\lim_{n \in \mathbb{N}_0} \bar{b}_{n,j} = b_j$ for all $j$, and

$$\lim_{n \in \mathbb{N}_0} \sum_{j=1}^{j_n} [\bar{b}_{n,j} - b_j]^2 = 0$$

for some subsequence $j_n \to \infty$. By thinning the subsequence $\mathbb{N}_0$, if necessary, we may suppose that $j_n$, $n \in \mathbb{N}_0$ are strictly increasing. There is a strictly decreasing sequence $1 > q_1 > q_2, \ldots$ for which $\lim_{n \in \mathbb{N}_0} n q_{j_n} = 0$. Let $p_j = q_j - q_{j+1}$ and let $F_j$ be the distribution which assigns mass $\frac{1}{2} p_j$ to $\pm 1/\sqrt{p_j}$ and mass $1 - p_j$ to $0$. With this choice of $F_1, F_2, \ldots$, let

$$\tilde{M}_{n,n} = \sum_{j=1}^{j_n} \bar{b}_{n,j} \zeta_{n,j}.$$

Then $P[\zeta_{n,j} \neq 0] \leq n p_j$, and

$$P[M_{n,n} \neq \tilde{M}_{n,n}] \leq n q_{j_n} \to 0$$

as $n \to \infty$ in $\mathbb{N}_0$. So, $\tilde{M}_{n,n}/\sqrt{n}$ has a limiting normal distribution with mean $0$ and variance $\kappa^2$ and, therefore,

$$\liminf_{n \in \mathbb{N}_0} \sum_{j=1}^{j_n} \bar{b}_{n,j}^2 = \liminf_{n \in \mathbb{N}_0} \frac{1}{n} E(\tilde{M}_{n,n}^2) \geq \kappa^2$$

and

$$\lim_{n \in \mathbb{N}_0} \sum_{j=j_n+1}^{\infty} \bar{b}_{n,j}^2 = 0.$$

It follows easily that $\lim_{n \in \mathbb{N}_0} \bar{\mathbf{b}}_n = \mathbf{b}$ in $\ell^2(\mathbb{N})$, and since $\mathbf{b}$ was an arbitrary weak limit point, this establishes the first assertion.

The second assertion is now immediate. Setting all of the variances but one to zero shows that $\lim_{n \to \infty} \bar{b}_{n,j}^2$ exists for a fixed $j$, in which case $\lim_{n \to \infty} \bar{b}_{n,j}$ exists, since $|b_{n+1,j} - \bar{b}_{n,j}| = O(\sqrt{n})$, as in the proof of Corollary 2. It then follows that $\bar{\mathbf{b}}_n$ converges weakly, from which the assertion follows since $\bar{\mathbf{b}}_n$, $n \geq 1$, are pre-compact. $\square$

**Acknowledgments.** It is a pleasure to acknowledge helpful discussions with Dalibor Volný. Examples like those in [14] and [18] provided useful insight.

## REFERENCES

[1] BORODIN, A. N. and IBRAGIMOV, I. A. (1994). Limit theorems for functionals of random walks. *Trudy Mat. Inst. Steklov.* **195**. [Translated into English: *Amer. Math. Soc.* (1995) Chapter 4, Sections 7–9]. MR1368394

[2] BROCKWELL, P. J. and DAVIS, R. A. (1991). *Time Series: Theory and Methods.* Springer, New York. MR1093459

[3] DEDECKER, J. and MERLEVÈDE, F. (2002). Necessary and sufficient conditions for the conditional central limit theorem. *Ann. Probab.* **30** 1044–1081. MR1920101

[4] DEDECKER, J., MERLEVÈDE, F. and VOLNÝ, D. (2007). On the weak invariance principle for non adapted sequences under projective criteria. *J. Theoret. Probab.* **20** 971–1004. MR2359065

[5] DENKER, M. and GORDIN, M. I. (1998). The central limit theorem for random perturbations of rotations. *Probab. Theory Related Fields* **111** 1–16. MR1626762

[6] DERRIENNIC, Y. and LIN, M. (1996). Sur le théorème limite central de Kipnis et Varadhan pour les chaînes réversibles ou normales. *C. R. Acad. Sci. Paris Sér. I Math.* **323** 1053–1057. MR1423219

[7] DERRIENNIC, Y. and LIN, M. (2001). The central limit theorem for Markov chains with normal transition operators, started at a point. *Probab. Theory Related Fields* **119** 508–528. MR1826405

[8] DUNFORD, N. and SCHWARTZ, J. (1958). *Linear Operators.* **I.** Interscience, New York. MR0117523

[9] DURIEU, O. and VOLNÝ, D. (2008). Comparison between criteria leading to the weak invariance principle. *Ann. Inst. H. Poincaré Probab. Statist.* **44** 324–340.

[10] GORDIN, M. I. (1969). The central limit theorem for stationary processes. *Soviet Math. Dokl.* **10** 1174–1176. MR0251785

[11] GORDIN, M. I. and HOLZMANN, H. (2004). The central limit theorem for stationary Markov chains under invariant splittings. *Stoch. Dyn.* **4** 15–30. MR2069365

[12] GORDIN, M. and LIFŠIC, B. A. (1981). A remark about a Markov process with normal transition operator. In *Third Vilnius Conference on Probability Theory and Mathematical Statistics* **1** 147–148. Akad. Nauk Litovsk. SSR, Vilnius. (In Russian.)

[13] KIPNIS, C. and VARADHAN, S. R. S. (1986). Central limit theorem for additive functionals of reversible Markov processes and applications to simple exclusions. *Comm. Math. Phys.* **104** 1–19. MR0834478

[14] KLICNAROVÁ, J. and VOLNÝ, D. (2007). Exactness of a Wu–Woodroofe's approximation with linear growth of variances. Submitted.

[15] MAXWELL, M. and WOODROOFE, M. (2000). Central limit theorems for additive functionals of Markov chains. *Ann. Probab.* **28** 713–724. MR1782272

[16] TONG, H. (1990). *Non-Linear Time Series: A Dynamical System Approach.* Oxford Univ. Press. MR1079320

[17] VOLNÝ, D. (1993). Approximating martingales and the central limit theorem for strictly stationary processes. *Stochastic Process. Appl.* **44** 41–74. MR1198662

[18] VOLNÝ, D. (2007). Exactness of martingale approximation and the central limit theorem. Submitted.

[19] WOODROOFE, M. (1992). A central limit theorem for functions of a Markov chain with applications to shifts. *Stochastic Process. Appl.* **41** 33–44. MR1162717

[20] WU, W. B. and WOODROOFE, M. (2004). Martingale approximations for sums of stationary processes. *Ann. Probab.* **32** 1674–1690. MR2060314

[21] ZHAO, O. and WOODROOFE, M. (2008). Law of the iterated logarithm for stationary processes. *Ann. Probab.* **36** 127–142. MR2370600

DIVISION OF BIOSTATISTICS
YALE UNIVERSITY
SCHOOL OF MEDICINE
NEW HAVEN, CONNECTICUT
USA

DEPARTMENT OF STATISTICS
UNIVERSITY OF MICHIGAN
462 WEST HALL
ANN ARBOR, MICHIGAN
USA